\documentclass[12pt]{amsart}
\usepackage{amssymb}

\begin{document}

\title
[Combinatorics of Topological Posets]
{Combinatorics of Topological Posets:\\
        Homotopy complementation formulas}

\author{Rade T. \v Zivaljevi\' c}
\address{\hskip-\parindent
        R. \v Zivaljevi\' c\\
        Mathematics Institute
        SANU\\
        Knez Mihailova 35/1, p.f. 367\\
        11001 Belgrade, Yugoslavia}
\email{rade@turing.mi.sanu.ac.yu}

\thanks{The main ideas of the paper were born during
the workshop ``Geometric Combinatorics''
(MSRI, February 1997). We would like to thank the organizers
for the support and the organization of this inspiring
conference. Research at MSRI is supported in part by NSF grant DMS-9701755.}

\begin{abstract}
We show that the well
known {\em homotopy complementation formula} of
Bj\" orner and Walker admits several closely related
generalizations on different classes of topological
posets (lattices). The utility of this technique is demonstrated on
some classes of topological posets including the Grassmannian
and configuration posets, $\widetilde{\mathcal  G}_n(R)$ and
${\rm exp}_n(X)$ which were introduced and studied by V.~Vassiliev.
Among other applications we present a reasonably complete description,
in terms of more standard spaces, of homology types of configuration
posets ${\rm exp}_n(S^m)$ which leads to a negative answer to a
question of Vassilev raised at the workshop ``Geometric
Combinatorics'' (MSRI, February 1997).
\end{abstract}

\maketitle

\newtheorem{exam}{Example}[section]
\newtheorem{theo}[exam]{Theorem}
\newtheorem{rem}[exam]{Remark}
\newtheorem{prob}[exam]{Problem}
\newtheorem{defin}[exam]{Definition}
\newtheorem{prop}[exam]{Proposition}
\newtheorem{cor}[exam]{Corollary}
\newtheorem{conj}[exam]{Conjecture}
\newtheorem{lema}[exam]{Lemma}


\section{Introduction}
\label{Intro}

One of the objectives of this paper is to initiate the study
of topological (continuous) posets and their order complexes
from the point of view
of Geometric Combinatorics. Recall that finite or more generally
locally finite partially ordered sets
(posets) already occupy one of privileged  positions in this field.
The well-known identity (Philip Hall)
$$
\mu(P) = \tilde\chi(\Delta(P))
$$
in spite of its simplicity, often serves as a good initial
example illustrating both the combinatorial and the
geometric (topological) nature of these objects.
As usual, $\mu(P)$ is the M\" obius number of the
poset $P$, $\Delta(P)$ is equally well known
{\em order complex} of all chains in $P$ and $\tilde\chi(K)$ is the
reduced Euler characteristic of the space $K$.
\par
The order complex construction
is of fundamental importance for geometric combinatorics.
Recall some of the highlights.
As demonstrated by the
Goresky-MacPher\-son formula and its subsequent
generalizations, \cite{GorPhe88, OrlTer, ZZTop},
the homology (homotopy) type of interesting spaces associated to an affine
arrangement ${\mathcal  A}$ in $R^n$ can be described in terms of
order complexes of the cones in the associated intersection lattice
$L({\mathcal  A})$. The problem of finding combinatorial formulas for
Pontryagin classes can be reduced,  \cite{GelMac92},
to the study of combinatorial models of Grassmannians which in turn are
defined as the order complexes of posets of classes of oriented matroids.
Order complexes, as combinatorial models of spaces associated to
loop spaces, appear in \cite{BKS91}.
The homotopy colimit construction,
viewed as a natural generalization of the order complex
construction, has been applied in \cite{WZZ, ZZTop} on diagrams of spaces
over posets  which led to new results about affine, projective and
other arrangements and to new insights about toric varieties, deleted
joins etc.
\par
Our final example brings us
closer to the immediate objectives of this paper.
Victor Vassiliev used systematically and skilfully his
"geometric resolution"  method, \cite{Vas91, Vas94, Vas97},
together with a variety of other techniques, to obtain far reaching
results in several mathematical fields. Some of these
applications, especially the theory of Vassiliev knot invariants
and invariants of ornaments, reveal that the geometric resolution
is in some cases very close to or generalizes
the order complex construction.
One of surprising Vassiliev's discoveries, \cite{Vas91, Vas97},
is that the order complex construction
$\Delta({\mathcal  P})$, properly interpreted, generalized
and applied to interesting {\bf topological posets} ${\mathcal  P}$
leads to interesting and elegant, geometric observations.
\par
In this paper we focus out attention on topological posets and their
order complexes, their homotopy and homology types
as combinatorial objects interesting in their own right.
We raise a general problem of determining
which combinatorial results about (discrete) posets have
interesting continuous analogs. The emphasis
is of course on results which belong to the field of
geometric combinatorics, more precisely the results
which do not necessarily have a
proper analogue (generalization) if the posets ${\mathcal  P}$ is
replaced by a more general topological category ${\mathcal  C}$.
\par
Our central result in this direction is that the well
known {\em homotopy complementation formula} of
Bj\" orner and Walker, \cite{Bjo91, BjoWal83},
which deals with the order complex of a finite lattice,
admits several closely related generalizations on
different classes of topological posets (lattices).
\par
These generalized formulas lead to explicit ``computations'' of
homotopy (homology) types of order
complexes of natural topological posets, including posets associated to
Grassmannians and configuration spaces $F(X,n)$.
As an illustration we show how some recent results of Vassiliev
mentioned above can be obtained by this method. Another application
is a reasonably complete description, in terms of
more standard spaces, of homology types of
configuration posets ${\rm exp}_n(S^m)$ of spheres, known
previously in the case of $S^1$. As a consequence we obtain
a negative answer to a question of Vassiliev, \cite{Vas97b},
who asked whether the order complex $\Delta({\rm exp}_n(X))$ is
homeomorphic to the join $X^{\ast (n)} = X\ast\ldots\ast X$.
\par
The paper is organized as follows. In section~\ref{main}
we give a brief review of some of the central results of the paper
emphasizing the link with the usual ``homotopy complementation formula'' of Bj\" orner and Walker.
In section~\ref{topposet} several different classes of topological
posets are introduced and analyzed. Most of these classes
naturally appeared in the course of the proof of homotopy
complementation formulas which occupy section~\ref{Homo}.
Applications to Grassmannian and configuration posets
and elementary proofs are presented in
section~\ref{Appl}.

\section{Main results}
\label{main}
Homotopy complementation formula of A.~Bj\" orner and
J.W.~Walker (\cite{BjoWal83}, see also theorem~\ref{Compl}) is a versatile tool for computing homotopy types of order complexes
$\Delta(P)$ of finite posets. For example let $\Pi_n$ be
the lattice of all partitions of the set $\{1,\ldots,n\}$,
ordered by the refinement relation, and $\widetilde\Pi_n$ the poset
obtained by deleting the maximum and the minimum elements,
$(1 2 \ldots n)$ and $(1)(2)\ldots(n)$. Assume $n\geq 3$.
Then the homotopy complementation formula implies
the following homotopy recurrence relation, \cite{BjoWal83}

\begin{equation}
\label{part_lat}
\Delta(\widetilde\Pi_n) \simeq \bigvee_{i=2}^n \Sigma(\Delta(\widetilde\Pi_{n-1}^i))
\end{equation}
\noindent
where $\widetilde\Pi_{n-1}^i$ is a lattice isomorphic to
$\widetilde\Pi_{n-1}$ and $\Sigma$ is the suspension operator.
>From here, it is easily deduced by induction that the
homotopy type of the
lattice $\Pi_n$, i.e. the homotopy type of the order complex
$\Delta(\widetilde\Pi_n)$, is the wedge of $(n-1)!$ copies of the sphere
$S^{n-3}$.
\par
Our first objective is to prove an analogous homotopy
complementation formula for topological posets
(lattices) which yields the following homotopy recurrence
relations. In these examples, $\widetilde{\mathcal  G}_n(R),
\widetilde{\mathcal  G}_n^{\pm}(R),
{\rm exp}_n(X)$ are the Grassmannian and configuration topological
posets defined in section~\ref{topposet} and
$\widetilde{\mathcal  B}_n := {\mathcal  B}_n\setminus\{\emptyset\}$ is the
the face poset of an $n$-simplex, i.e. the
(truncated) Boolean lattice on $\{0,1,\ldots ,n\}$.
\begin{equation}
\label{gras_lat}
\Delta(\widetilde{\mathcal  G}_n(R)) \simeq S^{n-1}\wedge
\Sigma(\Delta(\widetilde{\mathcal  G}_{n-1}(R)))
\end{equation}
\begin{equation}
\label{grasgras_lat}
\Delta(\widetilde{\mathcal  G}_n^{\pm}(R)) \simeq (S^{n-1}\vee S^{n-1})\wedge
\Sigma(\Delta(\widetilde{\mathcal  G}^{\pm}_{n-1}(R)))
\end{equation}
\begin{equation}
\label{conf_lat}
\Delta({\rm exp}_n(S^1)) \simeq S^{n}\wedge
(\Delta(\widetilde{\mathcal  B}_{n-1})/
\partial\Delta(\widetilde{\mathcal  B}_{n-1}))
\simeq S^{2n-1}
\end{equation}
These formulas should not be viewed as isolated examples, rather
they illustrate a general and simple proof scheme which is
potentially applicable in many different situations.
Very often, see section~\ref{Homo}, these recurrence
formulas have the form
\begin{equation}
\label{form}
\Delta({\mathcal  P}_n) \simeq P_n\wedge
\Sigma(\Delta({\mathcal  P}_{n-1}))
\end{equation}
where $P_n$ is a pointed ``parameter'' space. Specially  if
$P_n$ is a finite set we recover the usual wedge form
as in the example (\ref{part_lat}).
Examples (\ref{gras_lat}) and (\ref{conf_lat}) lead to new proofs
of results of Vassiliev, \cite{Vas91} and \cite{Vas97}, which was
an initial motivation for a general complementation
formula for topological posets.
The example (\ref{conf_lat}) is a special case of a general
formula of the form
\begin{equation}
\label{confconf_lat}
\Delta({\rm exp}_n(X)) \simeq {\rm Thom}_{n}(X\setminus\{x_0\})
\end{equation}
\noindent
where ${\rm Thom}_n(Y)$ is  the one-point compactification of
a vector bundle over a configuration space $B(Y,n)$ of $n$-element
subsets of $Y$.
Motivated by (\ref{conf_lat}), Vassiliev asked, \cite{Vas97b},
if an analogous formula
\begin{equation}
\Delta({\rm exp}_n(X)) \simeq X^{\ast(n)} = X\ast \ldots \ast X
\end{equation}
\noindent
holds for arbitrary topological spaces,
specially if the order complex of ${\rm exp}_n(S^2)$ is
homotopic to $S^{3n-1}$.
We show that the answer to this question is in general
negative and that formula \ref{confconf_lat} implies that
the conjecture ${\rm (6)}$ is false already in the case $X = S^2$.
More importantly, we are able to give a sufficiently
complete description of homology types of order complexes
$\Delta({\rm exp}_n(S^m))$ in terms of homologies
of some standard spaces.

\section{Topological posets}
\label{topposet}
\subsection{Motivating examples}
\label{examples}
A topological poset $({\mathcal  P},\leq, \tau)$ is a poset
$({\mathcal  P},\leq)$ and a Hausdorff topology $\tau$ on the
set ${\mathcal  P}$ such that the order relation $R_{\leq} :=
\{(p,q)\in {\mathcal  P}\times {\mathcal  P} \mid p\leq q \}$ is a closed subspace
of ${\mathcal  P}\times {\mathcal  P}$.
\par
A morphism $f : ({\mathcal  P}_1, \leq_1 ,\tau_1) \rightarrow
({\mathcal  P}_2, \leq_2, \tau_2)$ in the category
$TPos$ of topological posets is a continuous, monotone map
$f : {\mathcal  P}_1 \rightarrow {\mathcal  P}_2$.
\par
Before we introduce special classes of
topological posets and begin their analysis,
let us review some motivating examples.

\medskip

\begin{exam}
(Grassmannian posets, \cite{Vas91})
\label{grapo}
{\rm Suppose that $K$ is one of the classical (skew)
fields $R, C$ or $Q$.
The Grassmannian poset ${\mathcal  G}_n(K) = (G(K^n),\subseteq)$,
is the disjoint sum
$$ G(K^n) := \coprod_{i=0}^{n} G_i(K^n) $$
\noindent
where $G_i(K^n)$ is the manifold of all
$i$-dimensional linear subspaces of $K^n$.
The order in this poset is by inclusion,
$U\leq V$ iff $U\subseteq V$. Every Grassmannian poset
${\mathcal  G}_n(K)$ is a lattice
with the minimum element ${\hat 0} = \{0\}$, maximum
element ${\hat 1} = K^n$ and the rank function
$r : G(K^n) \rightarrow N$ defined by   $r(V) := {\rm dim}(V)$.
The poset $\widetilde{\mathcal  G}_n(K) := {\mathcal  G}_n(K)\setminus
\{\hat 0,\hat 1\}$ is called the truncated Grassmannian poset.
}
\end{exam}

\begin{exam}(Subspace posets)
{\rm
For a given space $X$, let
${\rm exp}(X)$ be the topological space of all closed
subspaces equipped with the Vietoris topology, \cite{Engel}.
The inclusion relation $\subseteq$ turns this space
into a topological poset. This poset is not very interesting itself,
at least from the point of view of combinatorics. Its major role is
to serve as a source of interesting subposets. For example circle,
polygon, disc etc. posets from \cite{AloSch88} are examples
of subposets of ${\rm exp}(R^d)$.
}
\end{exam}

\begin{exam}
\label{exam2}
(Configuration posets)
{\rm
Let ${\rm exp}_n(X) := \{A\in {\rm exp}(X) \mid \vert A\vert \leq n\}$
be the subposet of ${\rm exp}(X)$ consisting of all finite,
nonempty, closed subsets of $X$ of cardinality less or equal to
$n, \ n\in N$.  The space ${\rm exp}_n(X)$
was under the name $n$-th symmetric product of $X$ introduced by Borsuk
and Ulam, \cite{BoUl31}. Note that ${\rm exp}_n(X)$ is related to but not
the same as the orbit space $X^n/S_n$ where $S_n$ is the
symmetric group. The space ${\rm exp}_n(X)$
is viewed as a topological subposet of ${\rm exp}(X)$ and it
has a natural rank function $p : {\rm exp}_n(X)\rightarrow N$
defined by $p(A) = \vert A\vert$. Let $B(X,k) = p^{-1}(k)$ be
the space of all $k$-element subsets in ${\rm exp}_n(X)$.
If $F(X,k)$ is the usual configuration space  of all ordered
collections $(x_1,\ldots ,x_k)\in X^k, \, x_i\neq x_j$ for $i\neq j$,
then $B(X,k)$ is the space of all unordered, $k$-element
configurations, $B(X,k) \cong F(X,k)/S_k$.
}
\end{exam}

\begin{exam}
(Semialgebraic posets)
{\rm
Let $M$ be a semialgebraic set in $R^n$ and
$\leq$ an order relation on $M$ such that
$R_{\leq} = \{(x,y)\in M\times M \mid x\leq y \}$
is a semialgebraic subset of $M\times M$. These posets will be
generally called {\em semialgebraic posets} although in this generality
they often intersect with other classes of posets.
An example of a semialgebraic poset is $R^{n+1} \cong R^n\oplus R$
with the quadratic form $q(x,t) = x_1^2+\ldots + x_n^2 - t^2$
and the order relation defined by $(x,t)\leq (y,s)$
iff $q(x-y,t-s)\leq 0$. Other examples of semialgebraic posets
which have interesting applications in combinatorics can be found
in \cite{AloSch88}.
}
\end{exam}

\begin{exam}(Diagram posets)
\label{dop}
{\rm
Diagrams of spaces, specially diagrams
of spaces over finite posets entered combinatorics in
papers \cite{WZZ, ZZTop}. Recall that a diagram ${\mathcal  D}$
over a (discrete) poset $(P,\leq)$ is
a (contravariant) functor ${\mathcal  D}: P^{op}\rightarrow Top$ where $(P,\leq)$ is
viewed as a small category such that $p\rightarrow q$
iff $p\leq q$. In a more informal language a diagram
${\mathcal  D}$ consists of a family of spaces $\{D_p\}_{p\in P}$
and a family of ``connecting'' maps $\{d_{pq} : D_q\rightarrow
D_p\}_{p\leq q}$ such that $d_{pp} = 1_{D_p}$ and
$d_{pq}\circ d_{qr} = d_{pr}$ for $p\leq q\leq r$.
Every diagram ${\mathcal  D}$ gives rise to
a topological poset $(\widetilde{\mathcal  D},\preccurlyeq)$ on the space
$\widetilde{\mathcal  D} := \coprod_{p\in P}~{\mathcal  D}(p)$ where
for $x\in {\mathcal  D}(p)$ and $y\in {\mathcal  D}(q)$
$$
x\preccurlyeq y  \Longleftrightarrow (p\leq q) \quad {\rm and}\quad
f_{pq}(y) = x \, .
$$
}
\end{exam}

\subsection{The order complex of a topological poset}
\label{Order}
The order complex of a discrete poset $P$ is the simplicial complex
$\Delta(P)$ of all chains in $P$. If $({\mathcal  P},\leq,\tau)$, or simply
${\mathcal  P}$, is a topological poset then there exists
a space $\Delta({\mathcal  P})$ naturally associated
to ${\mathcal  P}$ which can serve as an `ordered complex'
of ${\mathcal  P}$. If the poset ${\mathcal  P}$ is interpreted as a small,
topological category $\widetilde {\mathcal  P}$ where
${\rm ob}(\widetilde {\mathcal  P}) = {\mathcal  P}$ and
${\rm mor}(\widetilde {\mathcal  P}) = \{(x,y)\in {\mathcal  P}^2 \mid
x\leq y\}$, then the space $\Delta({\mathcal  P})$ is naturally
homeomorphic to the classifying space
$B(\widetilde {\mathcal  P})$ of $\widetilde {\mathcal  P}$, \cite{Seg68}.
The reader who is not familiar with basics of this theory
may find it convenient to review first the definition~\ref{mirror2}
which gives an alternative and more elementary
way of defining the order complex
for a narrow but important class of topological $M$-posets.

\begin{defin}
\label{simplicial}
Given a topological poset ${\mathcal  P}$,
let $N_{\ast} = \{N_n({\mathcal  P})\}_{n=0}^{\infty}$
be the associated simplicial
space where $N_n({\mathcal  P})$ is the space of all (not necessarily
strictly) increasing, finite chains $x_0 \leq x_1 \leq\ldots\leq x_n$
in ${\mathcal  P}$, topologized as
a subspace of ${\mathcal  P}^{n+1}$.
The order complex $\Delta({\mathcal  P})$ of ${\mathcal  P}$ is
by definition the geometric realization,
\cite{Seg68, May72, May74, Seg74}, of the simplicial space
$N_{\ast}({\mathcal  P})$. More explicitly, the space $\Delta({\mathcal  P})$
is described as the union (colimit) of the inductively defined
sequence of spaces $F_p = F_p(\Delta({\mathcal  P}))$, where
$F_0 \cong {\mathcal  P}$ and  $F_p$ is constructed from
$F_{p-1}$ and $N_{\ast}$ by the following push-out diagram

$$
\begin{array}{ccc}
(N_p\times\dot\Delta^p)\cup(\delta N_p\times\Delta^p) &
\longrightarrow &  F_{p-1} \\
\downarrow & & \downarrow \\
N_p\times\Delta^p &\longrightarrow &  F_p
\end{array}
$$
The space $\delta N_p$ of all degenerated $p$-simplices
is in the case of posets
the space of all non-strictly increasing $p$-chains in ${\mathcal  P}$.

\end{defin}

Informally speaking, the space $\Delta({\mathcal  P})$
should be the union, as in the case of discrete posets,
of all simplices spanned by all finite,
{\em strictly increasing} chains in ${\mathcal  P}$.
This is indeed the case, since nondegenerated simplices
in the simplicial space $N_{\ast}({\mathcal  P})$
are in one--to--one
correspondence with strictly increasing finite chains in ${\mathcal  P}$.
In other words $\Delta({\mathcal  P})$ is isomorphic as a set,
with the usual order complex of ${\mathcal  P}$  seen as a discrete poset.
On the other hand, the topology on
${\mathcal  P}$ enters the definition of $\Delta({\mathcal  P})$
in an essential way and the following simple example should
make the difference perfectly clear.

\medskip
\noindent
{\bf Example}:
Let $(R^2,\leq)$ be the topological poset where $R^2$ has the
usual topology and $(x_1,y_1)\leq (x_2,y_2)$
iff $x_1\leq x_2$ and $y_1\leq y_2$.
Let $A = \{(x,y)\mid x+y=1, \ x,y \geq 0\}$
and let ${\mathcal  P} = \{(0,0)\}\cup A$
be a topological subposet of $R^2$. Then the topological
order complex $\Delta({\mathcal  P})$ is homeomorphic to the triangle
with vertices $(0,0), (1,0),(0,1)$ with the topology induced from
$R^2$ while the discrete order complex is the cone with vertex
$(0,0)$ over the (discrete) set $A$.

\medskip

The following elementary proposition,
connecting the order complex
of a {\em diagram poset} $\widetilde{\mathcal  D}$  (example~\ref{dop})
with the {\em homotopy colimit} (\cite{HolVog92, ZZTop, WZZ}) of
the diagram ${\mathcal  D}$,  serves as an initial
example for the naturality of the order complex construction.

\begin{prop}
If ${\mathcal  D} : Q\rightarrow Top$ is a diagram of spaces over a finite
poset $Q$ and $\widetilde{\mathcal  D}$ the associated topological poset
(example~\ref{dop}), then
$$
\Delta(\widetilde{\mathcal  D})\cong {\rm {\bf hocolim}}_Q\, {\mathcal  D}.
$$
\end{prop}

\subsection{Special classes of topological posets}

In this section we introduce some special classes of
topological posets. Our goal is to define classes
which are broad enough to include interesting posets
and sufficiently narrow to allow special technical arguments.
The first two classes, the classes of {\em A-posets} and
{\em B-posets}, are characterized as topological posets satisfying
some cofibration conditions. The other two classes of
{\em C-posets} and {\em M-posets} are much more special
but have many useful properties.
Their detailed analysis is given in section~\ref{MandC}.

\begin{defin}
\label{A-posets}
($A$-posets)
A topological poset ${\mathcal  P}$ is an $A$-poset if the associated
simplicial space $N_{\ast}({\mathcal  P})$
has the property that
$(N_n({\mathcal  P}), \delta N_n({\mathcal  P}))$ is a cofibration pair
for each $n\in N$ where  $\delta N_n({\mathcal  P})$ is the
space of degenerated $n$-simplices.
A simplicial space with this or other closely related property is often
in the literature referred to as {\em cofibrant}, {\em good} or {\em proper},
\cite{Seg68, May72, May74, Seg74}.
\end{defin}

\begin{defin}
\label{B-posets}
($B$-posets)
A topological poset ${\mathcal  P}$ is a $B$-poset or an {\em intervally
cofibrant} poset if for each two closed intervals ${\mathcal  I}_1,
{\mathcal  I}_2\subset {\mathcal  P}$, if ${\mathcal  I}_1\subset {\mathcal  I}_2$
then the inclusion $\Delta({\mathcal  I}_1)\hookrightarrow
\Delta({\mathcal  I}_2)$ is a cofibration. A closed interval in
${\mathcal  P}$ is by convention any set of the form
$[a,b]_{\mathcal  P}, {\mathcal  P}_{\geq a}, {\mathcal  P}_{\leq b}$
including the poset ${\mathcal  P}$ itself.
\end{defin}

\begin{rem}
\label{prim}
{\rm
The inclusion of $CW$-complexes is always a cofibration. It is
well known that a finite family ${\mathcal  F}$ of (semi)algebraic
sets in $R^n$ admits a compatible triangulation,
\cite{Hiro75, Loj64}. These two facts together are often sufficient,
especially in the context of $M$-posets defined bellow,
to conclude that a given topological poset is an $A$ or
$B$-poset. Alternatively for the same purpose one can use
well known general properties of cofibrations,
\cite{Lil73, Strom68, May72, Bred}.
}
\end{rem}

The class of $M$-posets
introduced in the following definition includes Grassmannian and Diagram
posets but the class of configuration posets is out of its scope.

\begin{defin}
\label{mirror1}
Let ${\mathcal  P}$ be a topological poset.
A pair $(Q,\mu)$ where $Q$ is a finite poset and $\mu
: {\mathcal  P}\rightarrow Q$ is a monotone map, is called a
{\em mirror} of ${\mathcal  P}$ if
\begin{enumerate}
\item
for all $p, q\in {\mathcal  P}$, if $p < q$ then $\mu(p) < \mu(q)$,
\item
for all $q\in Q, \, \mu^{-1}(q)$ is a nonempty closed subspace of ${\mathcal  P}$.
\end{enumerate}
In this case ${\mathcal  P}$ is called an $M$-poset over $Q$ and
$\mu : {\mathcal  P}\rightarrow Q$ is an associated mirror map or shortly
an $M$-map from ${\mathcal  P}$ to $Q$.
\end{defin}

By definition~\ref{mirror1}, every strictly increasing
chain $x_1 < \ldots < x_k$ in ${\mathcal  P}$ is by an $M$-map
sent to a strictly increasing chain in $Q$ of the same length.
This useful property of $M$-posets is crucial in the following definition
and the subsequent proposition~\ref{jojo}.

\begin{defin}
\label{mirror2}
Assume that ${\mathcal  P}$ is an $M$-poset with
$\mu : {\mathcal  P}\rightarrow Q$ as an associated mirror map.
Let ${\mathcal  P}^{\ast (Q)} := {\mathcal  J}_{i\in Q}~{\mathcal  P}_i$
be the join of the family $\{{\mathcal  P}_i\}_{i\in Q}$ of spaces
where ${\mathcal  P}_i \cong {\mathcal  P}$
for all $i$. Each element $a\in {\mathcal  P}^{\ast(Q)}$ has the form
$x = \Sigma_{i\in Q}~\lambda_ia_i$ where $a_i\in {\mathcal  P}_i$,
$\lambda_i\geq 0$ and $\Sigma_{i\in Q}\lambda_i = 1$.
Let ${\rm supp}(a) = \{i\in Q \mid \lambda_i > 0\}$.
Define the $\mu$-order complex  $\Delta_{\mu}({\mathcal  P})$
of ${\mathcal  P}$ as  the subspace of
${\mathcal  P}^{\ast (Q)}$ where
$$ a\in \Delta_{\mu}({\mathcal  P})
\Longleftrightarrow
{\rm supp}(a)\ \ \mbox{\rm\small is a chain and if} \
i < j \  \mbox{\rm\small in} \ {\rm supp}(a) \ \mbox{\rm\small
then} \ a_i < a_j
\  {\rm in} \ {\mathcal  P}
$$
\end{defin}

\begin{prop}
\label{jojo}  Suppose ${\mathcal  P}$ is an $M$-poset
and $\mu : {\mathcal  P}\rightarrow Q$ an associated $M$-map.
Then the spaces $\Delta_{\mu}({\mathcal  P})$ and $\Delta({\mathcal  P})$
are homeomorphic.
\end{prop}

\medskip
\noindent
{\bf Proof}: The proof is based on the fact that for each
$M$-poset ${\mathcal  P}$, the space
$N_n({\mathcal  P})\setminus \delta N_n({\mathcal  P})$ of
all nondegenerated simplices in ${\mathcal  P}$ is a closed subspace of
$N_n({\mathcal  P})$ for each $n\in N$. \hfill$\square$

\medskip
As already noted above, there are two situations when the definition of
$\Delta_{\mu}({\mathcal  P})$ seems to be specially convenient.

\begin{enumerate}
\item ${\mathcal  P}$ is a topological poset with a finite
rank function $\rho : {\mathcal  P}\rightarrow [n],
\ [n] =\{0,1,\ldots,n\}$ as a mirror function.
An example of such a poset is the Grassmannian poset
${\mathcal  G}_n(K) = (G(K^n),\subseteq)$.
\item  ${\mathcal  D} : Q\rightarrow Top$ is a diagram
of spaces. If $(\widetilde{\mathcal  D},\preccurlyeq)$ is the associated
topological poset then the obvious monotone
map $\mu : \widetilde{\mathcal  D}\rightarrow Q$ is a mirror map.
\end{enumerate}

\medskip
The first condition in the definition of a
mirror poset, definition~\ref{mirror1}, looks quite restrictive.
If this condition is deleted, we obtain a much broader
class of topological posets which still preserve some of the
favorable properties of $M$-posets.

\begin{defin}
A topological poset $({\mathcal  P},\leq)$ is called a $C$-poset,
or more precisely a $C$-poset over a finite poset $Q$,
if there exists a monotone map $\alpha : {\mathcal  P}\rightarrow Q$,
called a $C$-map, such that $\alpha^{-1}(q)$ is a nonempty closed subspace
of ${\mathcal  P}$ for each $q\in Q$.
\end{defin}

\subsection{Mirrors and diagrams}
\label{MandC}
The main result of this section, proposition~\ref{mirror3},
is that there exists a functor from the category of $M$-posets
over a fixed finite set $Q$ to the category of diagrams
over another poset, closely related to $Q$.
This is a useful observation because the replacement of a
topological poset (topological category) with a {\em non--discrete}
set of objects, by a diagram over
a {\em discrete} poset, allows us to use a variety of results
which are not available for general topological categories,
\cite{HolVog92, ZZTop, WZZ}.

\medskip
Recall (example~\ref{dop}) that every diagram of spaces
${\mathcal  D} : Q \rightarrow Top$ over a finite poset $Q$
yields a topological poset $\widetilde{\mathcal  D}$.
We ask when the converse is true.

\medskip
\noindent
{\bf Question}: Given a topological poset ${\mathcal  P}$ and a
mirror map $\mu : {\mathcal  P} \rightarrow Q$, how
far is ${\mathcal  P}$ from being a topological poset
$\widetilde{\mathcal  D}$ associated to a diagram
${\mathcal  D} : Q\rightarrow Top$ over $Q$.
More generally, when is it possible to define a diagram
${\mathcal  E} : R\rightarrow Top$ over a finite poset $R$,
constructed naturally from $Q$, such that the order complexes
$\Delta({\mathcal  P})$ and $\Delta(\widetilde{\mathcal  E})$
are naturally homeomorphic.

\medskip
Here are two examples that help us predict the answer.

\medskip
\noindent
{\bf Example}: Let $X = X_0\ast\ldots\ast X_n$ be the join
of a sequence $X_0,\ldots , X_n$ of spaces. Then
$X\cong \Delta({\mathcal  P})$ where ${\mathcal  P} := \coprod_{i=0}^n~X_i$
and for $x\in X_i$ and $y\in X_j$,
$x\lneqq y$ if and only if $i<j$. The obvious rank function
$r : {\mathcal  P}\rightarrow [n], \, [n] := \{0,1,\ldots ,n\}$,
is a mirror map. Let $R:=\Delta([n])$ be the order complex of
$[n]$. Then $R\cong {\mathcal  B}([n])\setminus\{\emptyset\}$
is the face poset of an $n$-simplex. Let ${\mathcal  E} : R^{op}
\rightarrow Top$ be the diagram defined by
${\mathcal  E}(A) := \prod_{i\in A}~X_i, \, A\in {\mathcal  B}([n])\setminus
\{\emptyset\}$, with obvious projections as the diagram maps. Then
$$
\Delta({\mathcal  P}) = {\bf hocolim}_Q {\mathcal  E} =
\Delta(\widetilde{\mathcal  E}).
$$

\begin{exam}
{\rm
Let ${\mathcal  P} = \widetilde{\mathcal  G}_n(R)$ be the truncated
Grassmannian poset defined in the example~\ref{grapo}.
Then the rank function $r : \widetilde{\mathcal  G}_n(R)
\rightarrow \langle n\rangle$ can be seen as a mirror map
from ${\mathcal  P}$ to the poset
$\langle n\rangle := \{1,\ldots ,n\}$. Let $R = \Delta(\langle n\rangle)$
be the order complex of $\langle n\rangle$ and assume that $R$ is ordered by inclusion.
Let ${\mathcal  E} : R^{op} \rightarrow Top$ be the diagram over $R$
defined by ${\mathcal  E}(I) := G_I(R^n)$ where $I = \{i_1<i_2<\ldots
<i_m\}$ is a nonempty subset (sequence) of $\langle n\rangle$
and $G_I(R^n) = G_{i_1,\ldots,i_m}(R^n)$ is the flag
manifold of all linear flags $F_1\subset\ldots\subset F_m$
in $R^n$ where ${\rm dim}(F_k) = i_k$ for all
$k=1,\ldots ,m$. The diagram map $e(J,I) : G_I(R^n)\rightarrow
G_J(R^n)$, for each pair of subsets $J\subset I\subset\langle n\rangle$,
is the obvious projection (restriction) map.
Then,
$$
\Delta(\widetilde{\mathcal  G}_n(R)) \cong {\rm {\bf hocolim}}_R\, {\mathcal  E}
\cong \Delta(\widetilde{\mathcal  E}).
$$
}
\end{exam}

\begin{defin}
\label{ediag}
Suppose that ${\mathcal  P}$ is an $M$-poset and $\mu :
{\mathcal  P}\rightarrow Q$ an associated mirror map.
Let $R = \Delta(Q)$ be the order complex of $Q$,
seen as a poset of faces of $\Delta(Q)$ ordered by
inclusion. There exists a naturally defined diagram of spaces
${\mathcal  E} : R^{op}\rightarrow Top$ over $R$ defined as follows.
Let $I = \{q_1<q_2<\ldots<q_m\}$ be a chain in $Q$ representing an
element of $R$. Define ${\mathcal  C}(I) := \prod_{k=1}^m~\mu^{-1}(q_k)$
and let ${\mathcal  E}(I)$ be the, possibly empty,
subspace of ${\mathcal  C}(I)$ defined by
$$
x = (x_1,\ldots, x_m) \in {\mathcal  E}(I) \quad \Longleftrightarrow
\quad x_1 < x_2 < \ldots < x_m \quad {\rm in} \quad {\mathcal  P}.
$$
For a subchain $J$ of $I$ let $e(J,I) : {\mathcal  E}(I)
\rightarrow {\mathcal  E}(J)$ be the projection (restriction) map
which restricts $x = (x_1,\ldots, x_m)$ to the subchain $J$.
\end{defin}

\begin{prop}
\label{mirror3}
Let ${\mathcal  P}$ be an $M$-poset with a mirror map
$\mu : {\mathcal  P}\rightarrow Q$. Let $R = \Delta(Q)$ be the
ordered complex of $Q$, viewed as the face poset
ordered by the inclusion and let
${\mathcal  E}: R^{op}\rightarrow Top$ be the diagram from
definition~\ref{ediag}. Then,
$$
\Delta({\mathcal  P})\cong {\rm {\bf hocolim}}_R\, {\mathcal  E} \cong
\Delta(\widetilde{\mathcal  E}).
$$
\end{prop}

\medskip\noindent
{\bf Proof}: From the assumption that $\mu^{-1}(q)$ is closed
for all $q\in Q$, we deduce that, in the spirit of
proposition~\ref{jojo}, the space $\Delta({\mathcal  P})$ can
be identified as a subspace of the join
${\mathcal  J}_{q\in Q}~\mu^{-1}(q)$. The rest of the proof
follows by inspection. \hfill$\square$

\medskip
As an application of proposition~\ref{mirror3}, we show how one
can get the information about the fundamental group
$\pi_1(\Delta({\mathcal  P}))$ of the order complex of an $M$-poset
${\mathcal  P}$ from the associated diagram ${\mathcal  E}$.

\begin{prop}
Let ${\mathcal  P}$ be an $M$-poset with a mirror map
$\mu : {\mathcal  P}\rightarrow Q$. Let $R := \Delta(Q)$
be the poset and ${\mathcal  E} : R^{op}\rightarrow Top$ the
associated diagram of spaces from definition~\ref{ediag}
and proposition~\ref{mirror3}. Then
$$
\pi_1(\Delta({\mathcal  P})) \cong {\bf colim}_R \ \pi_1({\mathcal  E})
$$
where $\pi_1({\mathcal  E}) : R^{op}\rightarrow Group$ is the diagram of
groups, $\pi_1({\mathcal  E})(I) :=
\pi_1({\mathcal  E}(I)),
\, I\in R$.
\end{prop}

\medskip
\noindent
{\bf Proof}: The result follows from proposition~\ref{mirror3}
and the Seifert-van Kampen theorem which for diagrams
over posets, cf. Brown \cite{Brown} p. 206, reads as follows
$$
\pi_1({\bf hocolim}_R \, {\mathcal  E}) \cong {\bf colim}_R \,
\pi_1({\mathcal  E}) .
$$

We continue this section with a discussion of $C$-posets,
their comparison with $M$-posets and some useful results which will be used
in section~\ref{Homo}.

\begin{prop}
\label{Avala}
Suppose that $({\mathcal  P},\leq)$ is a $C$-poset over $Q$ and let
$\alpha :{\mathcal  P}\rightarrow Q$ be an associated $C$-map.
Suppose that the order complex
$\Delta(D_q)$ is compact for each subposet
$D_q := \alpha^{-1}(q)\subset {\mathcal  P},\, q\in Q$.
Then there exists an $M$-poset $({\mathcal  P}_{\rm M},\leq_{\rm M})$,
where ${\mathcal  P}_{\rm M} := \coprod_{q\in Q} \, \Delta(D_q)$,
and an $M$-map $\mu : {\mathcal  P}_{\rm M}\rightarrow Q$
such that
$$
\Delta({\mathcal  P}_{\rm M})\cong \Delta({\mathcal  P}).
$$
\end{prop}

\medskip
\noindent
{\bf Proof}: The assumption about $({\mathcal  P},\leq)$
being a $C$-poset over $Q$ can be restated as the fact that
there exists a decomposition ${\mathcal  P} = \coprod_{q\in Q}\, D_q$
of ${\mathcal  P}$ into closed, convex subposets $D_q$ such that
for $x\in D_q, y\in D_{q'}$, if $x\leq y$ then $q\leq q'$.
The compactness assumption on $\Delta(D_q),\, q\in Q$,
implies that $D_q$ are compact subspaces of the compact space
${\mathcal  P}$. By the {\em saturation} $({\mathcal  E},\leq_e)$
of $({\mathcal  P},\leq)$, we mean a new order relation $\leq_e$ on the space
${\mathcal  E}:={\mathcal  P}=\coprod_{q\in Q}\, D_q$ where,
for given $a\in D_q$ and $b\in D_{q'}$,
$$
a\leq_e b \enskip \Longleftrightarrow \enskip
(q=q' \enskip {\rm and} \enskip a\leq b)
\enskip {\rm or} \enskip q<q'.
$$
Obviously, $\Delta({\mathcal  E})\cong {\mathcal  J}_{q\in Q}\, \Delta(D_q)$ and
$\Delta({\mathcal  P})$ can be identified as a closed
(compact) subspace of $\Delta({\mathcal  E})$.
Define $({\mathcal  E}_{\rm E},\leq_{\rm E})$ as the topological poset
with ${\mathcal  E}_{\rm E}:=\coprod_{q\in Q}\, \Delta(D_q)$ and
for $x\in \Delta(D_q), y\in \Delta(D_{q'}), \, x\leq_{\rm E} y$ iff $x=y$ or $q<q'$.
Obviously
$$\Delta({\mathcal  E}_{\rm E})\cong\Delta({\mathcal  E})\cong
{\mathcal  J}_{q\in Q}\, \Delta(D_q).
$$
Let $({\mathcal  P}_{\rm M},\leq_{\rm M})$ be a subposet of $({\mathcal  E}_{\rm E},
\leq_{\rm E})$ where ${\mathcal  P}_{\rm M}:={\mathcal  E}_{\rm E} =
\coprod_{q\in Q}\,\Delta(D_q)$ and the order relation is defined
as follows. Recall that for $x\in D_q,\, {\rm supp}(x)$
is the unique minimal chain $C$ in $D_q$ such that $x\in\Delta(C)
\subset\Delta(D_q)$. Then for $x\in D_q$ and $y\in D_{q'}$,
$$
x\leq_{\rm M} y \enskip \Longleftrightarrow \enskip
(q=q' \wedge x=y) \enskip\bigvee\enskip
(q<q' \wedge \ \mbox{\rm supp}(x) \cup  \mbox{\rm supp}(y) \
\mbox{\rm is a chain in}\ {\mathcal  P}).
$$
Both $({\mathcal  E}_{\rm E},\leq_{\rm E})$ and
$({\mathcal  P}_{\rm M},\leq_{\rm M})$
are $M$-posets over $Q$ where for example the
$M$-map $\mu : {\mathcal  P}_{\rm M}\rightarrow Q$ is determined
by $\mu(\Delta(D_q)) =\{q\}$. The desired isomorphism
$$\Delta({\mathcal  P}_{\rm M})\cong\Delta({\mathcal  P})$$
is now transparent since both spaces are identified
with the same subspace of
$$
\Delta({\mathcal  E}_{\rm E})\cong\Delta({\mathcal  E})\cong
{\mathcal  J}_{q\in Q}\, \Delta(D_q).
$$

\medskip
An important class of $C$-posets arises by the so called
``Grothendieck construction'',
\cite{Thomason, HolVog92}.

\begin{defin}
\label{Groth}
Let ${\mathcal  D} : Q^{op}\rightarrow   TPos$ be a diagram of topological
posets over a finite poset $Q$. In other words ${\mathcal  D}$
consists of a collection $\{(D_q,\leq_q)\}_{q\in Q}$
of topological posets indexed by $q\in Q$ and a collection
of continuous, monotone maps $\{d_{qq'}\}_{q\leq q'},
\ d_{qq'}: D_{q'}\rightarrow D_q$, such that $d_{qq} =
{\rm id}_{D_q}$ and $d_{qq'}\circ d_{q'q''} = d_{qq''}$
for $q\leq q'\leq q''$.
The Grothendieck construction applied on ${\mathcal  D}$ yields
a new topological poset $({\mathcal  C},\leq_c)$
where ${\mathcal  C}:= \coprod_{q\in Q}\ D_q\times\{q\}$ and $(x,q) \leq_c
(y,q')$ iff $x\leq_q d_{qq'}(y)$.
\end{defin}

The following proposition is well known in the case of discrete categories
(posets), \cite{Thomason}.

\begin{prop}
\label{Maja}
Suppose ${\mathcal  D} : Q^{op}\rightarrow TPos$ is a diagram of compact
topological posets over a finite poset $Q$. Let $({\mathcal  C},
\leq_c)$ be the topological poset obtained from ${\mathcal  D}$ by the
Grothendieck construction, definition~\ref{Groth}.
Let $X_{\mathcal  D} : Q^{op}\rightarrow Top$ be the diagram of spaces
defined by
$$
X_{\mathcal  D}(q) = X_q :=\Delta({\mathcal  D}(q))
\enskip {\rm and} \enskip x_{qq'} = \Delta(d_{qq'}) : X_{q'}\rightarrow
X_q.
$$
Then there is a homotopy equivalence
$$
{\bf hocolim}_Q\ X_{\mathcal  D} \simeq \Delta({\mathcal  C}).
$$
\end{prop}

\medskip
\noindent
{\bf Proof}: There is another diagram of spaces
$Y : Q^{op}\rightarrow Top$ over $Q$ defined by
$Y_q := \Delta({\mathcal  C}_{\geq q})$ where
${\mathcal  C}_{\geq q}:= \coprod_{r\geq q}\ D_r \subset {\mathcal  C}$.
The inclusion $Y_{q'}\hookrightarrow Y_q, \ q\leq q'$,
is a cofibration so by the projection lemma in \cite{WZZ},
$\Delta({\mathcal  C})\cong {\bf colim}_Q\ Y$.
There is a map of diagrams $\alpha : X_{\mathcal  D}\rightarrow Y$
where $\alpha_q : X_q\rightarrow Y_q$
is the inclusion $X_q = \Delta(D_q)\hookrightarrow
\Delta({\mathcal  C}_{\geq q}) = Y_q$.
The map $\alpha_q$ is a homotopy equivalence which is deduced
from proposition~2.1 in \cite{Seg68}, see also  section~3.3
of \cite{WZZ}.
By the homotopy lemma (\cite{WZZ})
which is better known as the May-Tornehave-Segal
theorem, \cite{Seg74, May74, HolVog92}, the map $\bar{\alpha} :
{\bf hocolim}_Q\ X_{\mathcal  D} \rightarrow {\bf hocolim}_Q\ Y$
is also a homotopy equivalence. Finally,
$$
\Delta({\mathcal  C})\cong {\bf colim}_Q\ Y\simeq {\bf hocolim}_Q\ X_{\mathcal  D}.
$$

\medskip
The proposition~\ref{S-contractible}
is needed in the proof of theorem~\ref{Begle}.
The proof of this proposition is based based on
proposition~\ref{toto} which, together with its
corollary~\ref{mizimo} may have some
independent interest.

\begin{prop}
\label{toto}
Suppose that ${\mathcal  E} : R^{op}\rightarrow Top$ is a diagram
of spaces over a finite poset $R$.
Let $\Sigma^n(X)$ be the $n$-fold iterated suspension of the
space $X$ and $\Sigma^n({\mathcal  E})$ the diagram over
$R$ defined by $\Sigma^n({\mathcal  E})(r) := \Sigma^n({\mathcal  E}(r))$
for each $r\in R$.
Let ${\mathcal  D} : \widetilde{\mathcal  B}_n^{op} \rightarrow Top$
be a diagram over the truncated Boolean lattice
$\widetilde{\mathcal  B}_n := {\mathcal  B}_n\setminus\{\emptyset\}$
on $[n] = \{0,1,\ldots ,n\}$ defined by
$$
{\mathcal  D}(I) := \left\{
\begin{array}{cl}
\Delta(R), & \mbox{\rm for} \, I\neq [n]\\
{\bf hocolim}_R \ {\mathcal  E}, & \mbox{\rm for} \, I = [n]
\end{array}
\right.
$$
where $d_{I,J} : {\mathcal  D}(J)\rightarrow {\mathcal  D}(I)$
is defined by $d_{I, J} = id_{\Delta(R)}$ for $J\neq [n]$
and $d_{I,[n]}$ is the canonical projection map.
Then,
$$
{\bf hocolim}_R\, \Sigma^{n}({\mathcal  E})
\cong {\bf hocolim}_{\widetilde{\mathcal  B}_n} \, {\mathcal  D} .
$$
\end{prop}

\medskip
\noindent
{\bf Proof}: Let ${\mathcal  F} : R\times \widetilde{\mathcal  B}_n
\rightarrow Top$ be the diagram over the product poset
defined by
$$
{\mathcal  F}(r,I) := \left\{
\begin{array}{cc}
{\mathcal  E}(r), & \mbox{\rm if} \, I = [n] \\
\ast, & \mbox{\rm otherwise}
\end{array}
\right.
$$
The Segal's homotopy push down construction,
\cite{Seg74, WZZ, HolVog92}, applied on ${\mathcal  F}$
and two projections
$\pi_1 : R\times \widetilde{\mathcal  B}_n
\rightarrow R$ and $\pi_2 : R\times \widetilde{\mathcal  B}_n
\rightarrow \widetilde{\mathcal  B}_n$, yields two diagrams
${\mathcal  F}_1 :=(\pi_1)_{\ast}({\mathcal  F})$ and
${\mathcal  F}_2 := (\pi_2)_{\ast}({\mathcal  F})$,
over $R$ and $\widetilde{\mathcal  B}_n$.
In light of the homotopy pushdown theorem, \cite{HolVog92, WZZ},
it is sufficient to show that these two diagrams are locally
homotopy equivalent to diagrams $\Sigma^n({\mathcal  E})$
and ${\mathcal  D}$ respectively. More precisely,
we will show that there exist maps of diagrams
$\alpha : {\mathcal  F}_1 \rightarrow \Sigma^n({\mathcal  E})$
and $\beta : {\mathcal  F}_2 \rightarrow {\mathcal  D}$ such that
the maps $\alpha(r) : {\mathcal  F}_1(r)\rightarrow
\Sigma^n({\mathcal  E})(r)$ and $\beta(I) : {\mathcal  F}_2(I)
\rightarrow {\mathcal  D}(I)$ are homotopy equivalences.
By definition
$$
{\mathcal  F}_1(r) = {\bf hocolim}~{\mathcal  F}\vert\pi^{-1}(R_{\geq r})
$$
On the other hand the restriction diagram
${\mathcal  F}\vert (\{r\}\times \widetilde{\mathcal  B}_n)$
is a ``retract'' of  ${\mathcal  F}\vert {\pi^{-1}(R_{\geq r})}$.
The general principle, \cite{Seg74, HolVog92}
and section~3.3 of \cite{WZZ}, about homotopies
arising from natural transformations, permits us to conclude that
$$
\Sigma^n{\mathcal  E}(r) \cong {\bf hocolim}~{\mathcal  F}\vert (\{r\}\times
\widetilde{\mathcal  B}_n) \simeq {\mathcal  F}_1(r) =
(\pi_1)_{\ast}({\mathcal  F})(r)
$$
Similarly,
$$
{\mathcal  F}_2([n]) = (\pi_2)_{\ast}({\mathcal  F})([n])
= {\bf hocolim}~{\mathcal  F}\vert\pi_2^{-1}([n])
\cong {\bf hocolim}_R~{\mathcal  E}
$$
and for $I\neq [n], \, I\in  \widetilde{\mathcal  B}_n$,
$$
{\mathcal  F}_2(I) = (\pi_2)_\ast({\mathcal  F})(I)
\simeq {\bf hocolim}~{\mathcal  F}\vert(R\times\{I\})\cong \Delta(R)
$$
which completes the proof.
\hfill$\square$

\medskip
The following corollary, which is needed in the proof of
proposition~\ref{S-contractible}, follows from
the proposition~\ref{toto} by a single application of the
{\em homotopy} lemma, \cite{WZZ}.

\begin{cor}
\label{mizimo}
If the order complex $\Delta(R)$ is contractible than
the operations $\Sigma^n$ and ${\bf hocolim}$ commute
up to homotopy,
$$
{\bf hocolim}_R \, \Sigma^n{\mathcal  E} \simeq
\Sigma^n({\bf hocolim}_R \, {\mathcal  E}) .
$$
\end{cor}

\begin{prop}
\label{S-contractible}
Suppose ${\mathcal  D}: Q^{op}\rightarrow TPos$ is a diagram of compact topological
posets over the linearly  ordered set
$Q = \langle n\rangle := \{1,\ldots ,n\}$.
Suppose $D_i := {\mathcal  D}(i)$ is stably or $\Sigma$-contractible for each
$i = 1,\ldots ,n$ which means that the iterated suspension $\Sigma^{k}(\Delta(D_i))$ is
contractible for some $k\geq 0$.
Let $({\mathcal  C},\leq_c)$ be the poset from the
definition~\ref{Groth} and $\alpha : {\mathcal  C}\rightarrow Q$
the monotone map which turns ${\mathcal  C}$ into a $C$-poset over $Q$.
Then each fiber of the map
$\bar\alpha = \Delta(\alpha) :
\Delta({\mathcal  C})\rightarrow\Delta(\langle n\rangle)$
is $\Sigma$-contractible, specially cohomologically trivial.
\end{prop}

\medskip
\noindent
{\bf Proof}: The lemma is proved by induction on $n\in N$.
Note that $\Delta({\mathcal  C})$ is $\Sigma$-contractible
by propositions \ref{Maja} and \ref{toto}.
Assume, as the induction hypothesis,
that the fiber $\bar\alpha^{-1}(x)$ is $\Sigma$-contractible
for each element $x\in \partial\Delta(\langle n\rangle)$,
where  $\partial\Delta(\langle n\rangle)$ is the boundary of the simplex
$\Delta(\langle n\rangle) \cong \Delta^{n-1}$.
Suppose that $x$ belongs to the interior of $\Delta(\langle n\rangle)$.
Let ${\mathcal  C}_M$ be the $M$-poset associated to the $C$-poset
${\mathcal  C}$ constructed in the proposition~\ref{Avala}.
The poset $\Delta(\langle n\rangle)$, ordered by the reversed
inclusion, is in different notation the truncated
Boolean lattice $\widetilde{\mathcal  B}_{n-1} =
{\mathcal  B}_{n-1}\setminus\{\emptyset\}$ i.e. the face poset of an
$(n-1)$-dimensional simplex.
Then, cf. definition~\ref{ediag} and proposition~\ref{mirror3}, there exists a diagram of spaces
${\mathcal  E} : \widetilde{\mathcal  B}_n\rightarrow Top$
associated to the $M$-poset ${\mathcal  C}$ such that
${\bf hocolim}\ {\mathcal  E} \cong \Delta({\mathcal  C})$ and such that
${\mathcal  E}(I),\ I=\{i_1<\ldots <i_k\}$, can be identified
as the fiber $F_I$ of the map $\bar\alpha =\Delta(\alpha)$ over
the barycenter of the face $I$.
By applying the $\Sigma^m$ operation on the diagram
${\mathcal  E}$, where $m$ is big enough, and using the corollary~\ref{mizimo},
we obtain
$$
\ast \simeq {\bf hocolim}_{\widetilde{\mathcal  B}_n} \, \Sigma^m({\mathcal  E})
\simeq S^{n-1} \ast \Sigma^m(F_{\langle n\rangle})
\cong \Sigma^{m+n}(F_{\langle n\rangle})
$$
which means that the fiber $F_{\langle n\rangle}$ is
indeed $\Sigma$-contractible. \hfill$\square$

\section {Homotopy complementation formulas}
\label{Homo}

\begin{theo}{\rm (Homotopy Complementation Theorem,
Bj\" orner and Walker 1983)}
\label{Compl}
\par\noindent
Let $L$ be a bounded lattice and $z\in \widetilde L:= L - \{\hat 0,\hat 1\}$.
Let ${\rm Co}(z) := \{x\in \widetilde L \mid x\wedge z = \hat 0,
x\vee z = \hat 1\}$. Then
\begin{enumerate}
\item
The poset $\widetilde L\setminus {\rm Co}(z)$ is contractible
i.e. the order complex $\Delta(\widetilde L\setminus {\rm Co}(z))$
is contractible.
\item If ${\rm Co}(z)$ is an antichain, then
\begin{equation}
\label{wedge}
\Delta(\widetilde L) \simeq
\bigvee_{y\in {\rm Co}(z)} \Sigma(\Delta(\widetilde L_{<y})
\ast \Delta(\widetilde L_{>y})).
\end{equation}
\end{enumerate}
\end{theo}

The proof of this theorem consists of two parts. In the first part
it is shown that the poset $\widetilde L_1 := \widetilde L\setminus {\rm Co}(z)$
is contractible which implies that $\Delta(\widetilde L)\simeq
\Delta(\widetilde L)/\Delta(\widetilde L_1)$. In the second part
it is shown that the homotopy type of any space of the form
$\Delta(P)/\Delta(P\setminus C)$ where $P$ is a finite poset and
$C$ an antichain, specially in the case
$P = \widetilde L, C = {\rm Co}(z)$, has the desired wedge decomposition.
\par
Our objective is to extend this theorem to the case of
topological posets.
This goal is achieved in three steps, each step consisting of an
appropriate answer, or several answers to one of the following questions.

\begin{description}
\item[Q1.] Let ${\mathcal  P}$ be topological poset and
$C\subset {\mathcal  P}$ its subset, specially an
antichain. When is ${\mathcal  P}\setminus C$ contractible?
\item[Q2.]
If ${\mathcal  R}:= {\mathcal  P}\setminus C$ is contractible, when can we conclude
that $\Delta({\mathcal  P})$ and the quotient space
$\Delta({\mathcal  P})/\Delta({\mathcal  R})$ have the same homotopy type?
\item[Q3.]
Is there a decomposition formula, analogous to (\ref{wedge}),
for the homotopy
type of the space $\Delta({\mathcal  P})/\Delta({\mathcal  R})$?
\end{description}

\subsection{The first question}
\label{1and1}

\begin{prop}
\label{4-1}
Suppose that $({\mathcal  R}, \leq )$ is a topological poset
and let $z\in {\mathcal  R}$ such that the least upper bound $x\vee z$ exists
for each $x\in {\mathcal  R}$. Assume that the map
$\phi_z : {\mathcal  R}\rightarrow {\mathcal  R}, \, \phi_z(x) := x\vee z$
is continuous. Then the order complex $\Delta({\mathcal  R})$ of ${\mathcal  R}$
is contractible.
\end{prop}

\medskip
\noindent
{\bf Proof}: By proposition~2.1 from \cite{Seg68}, if
$f, g : {\mathcal  R}\rightarrow {\mathcal  R}$ a two morphisms in $TPos$
such that $f(x)\leq g(x)$ for each $x\in {\mathcal  R}$, then
the induced maps $\Delta(f), \Delta(g) : \Delta({\mathcal  R})
\rightarrow \Delta({\mathcal  R})$ are homotopic.
It follows that all three maps $id, \phi_z$ and
$c : {\mathcal  R}\rightarrow {\mathcal  R}, \, c(x) := z$,
induce homotopic maps of order complexes,
hence $\Delta({\mathcal  R})$ must be contractible.
\hfill$\square$

\medskip
\noindent
\begin{rem}
\label{4-2}
{\rm
The homotopy $H : \Delta({\mathcal  R})\times I\rightarrow \Delta({\mathcal  R})$
of two maps $\Delta(f), \Delta(g) : \Delta({\mathcal  R})\rightarrow
\Delta({\mathcal  R})$, where $f(x)\leq g(x)$ for each $x\in {\mathcal  R}$,
is defined as follows. Define
$F : {\mathcal  R}\times \{0<1\}\rightarrow {\mathcal  R}$ by
$F(x,0):= f(x)$ and $F(x,1):=g(x)$. Then,
$$
\Delta({\mathcal  R}\times\{0<1\}) \cong \Delta({\mathcal  R})\times
\Delta(\{0<1\}) \cong \Delta({\mathcal  R})\times I
$$
and $H := \Delta(F) : \Delta({\mathcal  R})\times I\rightarrow
\Delta({\mathcal  R})$.
More explicit description of $H$, at least in the case of
$M$-posets (definition~\ref{mirror1}) is
following. Let $C = \{x_0< \ldots < x_p\}$ be a chain
in ${\mathcal  R}$ and let $x\in \Delta(C)\subset\Delta({\mathcal  R})$
which, definition~\ref{mirror2}, has the form
$x = \Sigma_i \ \lambda_i x_i$.
Then $H$ can be described as a ``linear'' homotopy
between $\Delta(f)$ and $\Delta(g)$ where
$$
H(x,t) = (1-t)\Delta(f)(x) + t \Delta(g)(x) =
(1-t)\Sigma_i \ \lambda_i f(x_i) + t\Sigma_i \ \lambda_i g(x_i).
$$
}
\end{rem}

The continuity of the map $\phi_z : {\mathcal  R}\rightarrow {\mathcal  R}$
was used in essential way in the proof of proposition~\ref{4-1}.
In some cases it is still possible to reach the desired conclusion
assuming that $\phi_z$ is only partially continuous.
We illustrate the relevant ideas in the special case of the Grassmannian poset ${\mathcal  G}_n(R)$.

\begin{prop}
\label{4-3}
Let $\widetilde{\mathcal  G}_n(R)$ be the truncated Grassmannian poset,
example~\ref{grapo}. Given $z\in G_1(R^n)$, let ${\rm Co}(z)
:= \{l\in G_{n-1}(R^n) \mid z\not\subseteq l\}$ be the space of all
complements of $z$ in ${\mathcal  G}_n(R)$.
Define ${\mathcal  R} := \widetilde{\mathcal  G}_n(R) \setminus {\rm Co}(z)$.
Then the order complex $\Delta({\mathcal  R})$ is contractible.
\end{prop}

\noindent
{\bf Proof}: The map $\phi_z : {\mathcal  R}\rightarrow {\mathcal  R}, \
\phi_z(x) := x\vee z$, is no longer continuous so
it has to be modified.    Let ${\mathcal  Z}:= \{l\in {\mathcal  R}
\mid z\subseteq l\}$.
Let us observe first that the inclusion
$\Delta({\mathcal  Z})\hookrightarrow \Delta({\mathcal  R})$ is
a cofibration, cf. remark~\ref{prim}. We deduce (\cite{Bred})
that there exists a continuous function
$\alpha : \Delta({\mathcal  R})\rightarrow I$ such that
$\alpha^{-1}(0) = \Delta({\mathcal  Z})$ and
${\mathcal  U} := \alpha^{-1}([0,1))$ deforms to $\Delta({\mathcal  Z})$
through $\Delta({\mathcal  R})$, with $\Delta({\mathcal  Z})$ fixed.
The function $\phi_z$ is modified with the aid of the function
$\alpha$. Using the notation of the remark~\ref{4-2},
let $\psi : \Delta({\mathcal  R})\rightarrow
\Delta({\mathcal  R})$ be defined as follows.
For $x_i\in C$, let $\psi(x_i) := \alpha(x_i)(x_i\vee z) +
[1-\alpha(x_i)]x_i$. For a general $x\in\Delta(C), \, \psi$
is defined by linear extension,
$$
\psi(x) = \psi(\Sigma_i \lambda_i x_i) =
\Sigma_i \lambda_i \psi(x_i) .
$$
The continuity of this map follows from the fact that
$\phi_z$ is continuous away from ${\mathcal  Z}$ and the modification
takes care of points close or inside ${\mathcal  Z}$.
Note that $\psi$ does not arise as a map of the form
$\Delta(h)$ for some $TPos$-morphism $h : {\mathcal  R}\rightarrow {\mathcal  R}$.
On the other hand $\psi(x)$ is always on the ``interval''
connecting $\Delta(\phi_z)(x)$ and $x\in \Delta({\mathcal  R})$, hence
a continuous ``linear'' homotopy $H : \Delta({\mathcal  R})\times I
\rightarrow \Delta({\mathcal  R})$ is still well defined and
continuous. Note that $\psi$ maps $\Delta({\mathcal  R})$ into
the neighborhood ${\mathcal  U}$ which itself can be deformed into
$\Delta({\mathcal  Z})$. Finally, since the last space is contractible,
$\psi$ is null-homotopic, so is the identity map and
$\Delta({\mathcal  R})$ is contractible. \hfill$\square$

\medskip
In the following proposition we record for the future reference,
a more general statement which is proved along the lines
of the proof of proposition~\ref{4-3}.

\begin{prop} Suppose that ${\mathcal  R}$ is a compact topological poset
and let $z\in {\mathcal  R}$ be a minimal (maximal) element
such that the least upper bound $\phi_z(x) := x\vee z$
(the least lower bound $x\wedge z$) exists for each $x\in {\mathcal  R}$.
Suppose that the function $\phi_z : {\mathcal  R}\setminus
{\mathcal  R}_{\geq z}\rightarrow {\mathcal  R}$ is continuous and assume
that the inclusion map $\Delta({\mathcal  R}_{\geq z})\rightarrow
\Delta({\mathcal  R})$  is a cofibration. In the dual case
the upper cone ${\mathcal  R}_{\geq z}$ is replaced by the upper cone
${\mathcal  R}_{\leq z}$. Then the order complex
$\Delta({\mathcal  R})$ is contractible.
\end{prop}

\medskip
The following theorem is in some sense our most complete
answer to the first  question ({\bf Q1}). Its proof reveals that
this question is reducible to a form of the
Quillen's {\em fiber theorem}, \cite{Qui78, Bjo89-1, WZZ},
for topological posets. Quillen's result is known to hold,
\cite{HolVog92}, for diagrams of spaces over topological
categories with the discrete set of objects.
Unfortunately it doesn't seem to be available yet in the case
of general topological categories, specially in the case
of topological posets. This is the reason why we use instead
the Vietoris-Begle mapping theorem, \cite{Spanier}.

\begin{theo}
\label{Begle}
Let $({\mathcal  P},\leq)$ be a compact topological $B$ and $M$-poset,
definitions~\ref{B-posets} and \ref{mirror1}, such that the poset
$\hat{\mathcal  P} := {\mathcal  P}\cup\{\hat 0,\hat 1\}$ with added
minimum and maximum elements $\hat 0$ and $\hat 1$ is a lattice.
Let $z\in {\mathcal  P}$ and let ${\rm Co}(z)$ be the set of all
complements of $z$ in $\hat{\mathcal  P}$,
${\rm Co}(z) := \{x\in {\mathcal  P}\mid  x\wedge z = \hat 0 \
{\rm and} \  x\vee z = \hat 1\}$. Then the order complex
$\Delta({\mathcal  R})$ of the poset ${\mathcal  R}:= {\mathcal  P}\setminus
{\rm Co}(z)$ is cohomologically trivial,
i.e. has the cohomology of a point.
More generally, for the conclusion of the
theorem it is not necessary to assume that $\hat{\mathcal  P}$ is a
lattice. It suffices to assume that
$z\vee x$ and $z\wedge x$ exist in $\hat{\mathcal  P}$ for all $x\in {\mathcal  P}$.
\end{theo}
\noindent
{\bf Proof:} Let $e : {\mathcal  R}\rightarrow {\mathcal  R},
e(x) = x$, be the identity map and $c : {\mathcal  R}\rightarrow {\mathcal  R},
c(x) = z$ a constant map.
In order to show that $\Delta({\mathcal  R})$ is cohomologically trivial,
it suffices to prove that the maps ${\bar e} = \Delta(e)$ and
${\bar c} = \Delta(c)$ from $\Delta({\mathcal  R})$ to
$\Delta({\mathcal  R})$ induce the same homomorphism
$e^{\ast} = c^{\ast} : {\check H}^{\ast}(\Delta({\mathcal  R}))
\rightarrow {\check H}^{\ast}(\Delta({\mathcal  R}))$ of the associated
Alexander cohomology groups.
\par
Let ${\mathcal  R}\times {\mathcal  R}$ be the product poset and let
${\mathcal  F}:= \{(x,y)\in {\mathcal  R}\times{\mathcal  R}\mid
x\in {\mathcal  R} \ {\rm and} \ (x\leq y \ {\rm or} \ y\leq z)\}$.
Let $\pi_1, \pi_2 : {\mathcal  F}\rightarrow {\mathcal  R}$
be the obvious projection maps, $\pi_1(x,y)=x$ and $\pi_2(x,y)=y$,
restricted on ${\mathcal  F}$.
Define $d, j : {\mathcal  R}\rightarrow {\mathcal  F}$
by $d(x) = (x,x)$ and $j(x) = (x,z), x\in {\mathcal  R}$.
All maps $\pi_1, \pi_2, d, j$ are monotone and continuous
so they induce the corresponding maps, ${\bar \pi_1},
{\bar \pi_2}, {\bar j}, {\bar d}$
of ordered complexes. Obviously,
\begin{equation}
\label{comp}
\pi_1\circ d = \pi_1\circ j = e,  \,
e = \pi_2\circ d,
\,  c = \pi_2\circ j .
\end{equation}
The equality $e^{\ast} = c^{\ast} : {\check H}^{\ast}(\Delta({\mathcal  R}))
\rightarrow {\check H}^{\ast}(\Delta({\mathcal  R}))$
follows immediately if we prove that the homomorphism
$\pi^{\ast}_1 : {\check H}^{\ast}(\Delta({\mathcal  R}))
\rightarrow {\check H}^{\ast}(\Delta({\mathcal  F}))$,
associated to the projection $\pi_1 : {\mathcal  F}\rightarrow {\mathcal  R}$,
is an isomorphism.
Indeed, $d^{\ast} = j^{\ast}$ follows from $d^{\ast}\circ \pi^{\ast}_1
= j^{\ast}\circ \pi^{\ast}_1$ and $e^{\ast} = d^{\ast}\circ\pi^{\ast}_2
= j^{\ast}\circ\pi^{\ast}_2 = c^{\ast}$.
We shall demonstrate that $\pi^{\ast}_1$ is an isomorphism by showing that
${\bar\pi_1} : \Delta({\mathcal  F})\rightarrow \Delta({\mathcal  R})$
satisfies all conditions of the Vietoris-Begle mapping theorem,
\cite{Spanier}. The complex $\Delta({\mathcal  P})$ is compact
since ${\mathcal  P}$ is compact and ${\mathcal  P}$ is an $M$-poset.
Hence, the map $\bar\pi_1$ is closed as a continuous map
of compact spaces. It is surjective since the map
$\pi_1 : {\mathcal  F}\rightarrow {\mathcal  R}$ is surjective.
It remains to be proved that the fiber $\bar\pi^{-1}_1(w)$ is
cohomologically trivial for each $w\in\Delta({\mathcal  R})$.
We actually prove that the fiber $\bar\pi^{-1}(w)$ is
$\Sigma$-contractible (cf. proposition~\ref{S-contractible}).
\par
Let $w = \lambda_1 x_1 +\ldots +\lambda_n x_n$ where
$x_1 < x_2 <\ldots < x_n$ is a strictly increasing chain $C$
in ${\mathcal  R}$ and $\lambda_i > 0$ for all $i=1,\ldots ,n$.
Let ${\mathcal  E}:= \pi^{-1}_1(C)$ be the subposet of
${\mathcal  F}$ ``over'' the chain $C$. Let us show that
${\mathcal  E}$, together with the
monotone map $\bar\pi_1 : {\mathcal  E}\rightarrow C$, satisfies conditions
of proposition~\ref{S-contractible}.
Let $x\in\{x_1,\ldots ,x_n\}$. We are supposed to show that
$B(x) = \bar\pi^{-1}_1(x) = \{(x,t)\mid x\leq t \ {\rm or} \ t\leq z\}$
is contractible.
By assumption either $z\wedge x$ or $z\vee x$ exists in ${\mathcal  R}$.
Suppose $a:=z\wedge x$ exists, the other case
is treated analogously. Let $U(x) := \{t\in B(x)\mid
t\leq a \ {\rm or} \ a\leq t\}, \, V(x) := \{t\in B(x)\mid
t\leq z\}$ and $W(x) := U(x)\cap V(x)$.
Obviously $\Delta(B(x)) = \Delta(U(x)) \cup \Delta(V(x))$ and
$\Delta(U(x)) \cap \Delta(V(x)) = \Delta(W(x))$.
Here we use in an essential way the fact that $a$ is the greatest
lower bound of $x$ and $z$ so every chain in $B(x)$ which
intersects $U(x)\setminus V(x)$ is contained in $U(x)$.
All spaces $\Delta(U(x)),  \Delta(V(x)), \Delta(W(x))$
are contractible and the inclusion maps between them
are cofibrations since ${\mathcal  P}$ is a $B$-poset.
By the gluing lemma, \cite{Brown, WZZ},
$B(x)$ is also contractible.
Finally, an application of proposition~\ref{S-contractible} yields that
the fiber $\bar\pi^{-1}_1(a)$ is $\Sigma$-contractible.
This proves that $\pi^{\ast}_1$ is an isomorphism and the theorem follows.
\hfill$\square$

\subsection{The second and the third question}
\label{2and3}

Assuming that $\Delta({\mathcal  R})$ is contractible,
the spaces $\Delta({\mathcal  P})$ and
$\Delta({\mathcal  P})/\Delta({\mathcal  R})$
have the same homotopy type if the inclusion map
$\Delta({\mathcal  R}) \hookrightarrow \Delta({\mathcal  P})$
is a cofibration. Some of the methods how this condition
can be checked are discussed in the remark~\ref{prim}.
A general impression is that the answer to the question
{\bf Q2} is positive in all cases of interest.
\par
We focus our attention now to the third question.
Assume that the following assumptions hold
until the end of this section.

\begin{description}
\item[A1.]
$\Delta({\mathcal  P})$ is a compact,
\item[A2.]
${\mathcal  R}:={\mathcal  P}\setminus C$ where $C\subset {\mathcal  P}$ is
an open antichain.
\end{description}
It follows from {\bf A2} that $\Delta({\mathcal  R})$ is a closed
subspace of $\Delta({\mathcal  P})$. By {\bf A1}, the space
$Y := \Delta({\mathcal  P})\setminus\Delta({\mathcal  R})$ is
locally compact so the question {\bf Q3} reduces to the
analysis of the one-point compactification $\bar{Y}
:= Y\cup\{\infty\}$ of $Y$. Let us define auxiliary posets
$$
E_{{\mathcal  P},C} := \{(c,x)\in {\mathcal  P}^2\mid
c\in C \ \mbox{\rm and} \ (x\geq c \ \mbox{\rm or} \ x\leq c)\}
$$
$$
\dot E_{{\mathcal  P},C} := \{(c,x)\in {\mathcal  P}^2\mid
c\in C \ \mbox{\rm and} \ (x > c \ \mbox{\rm or} \ x < c)\}.
$$
The order complexes ${\mathcal  E}_{{\mathcal  P},C} := \Delta(E_{{\mathcal  P},C})$
and $\dot {\mathcal  E}_{{\mathcal  P},C} :=
\Delta(\dot E_{{\mathcal  P},C})$ of these two posets
can be naively seen as the disc and the sphere bundle
associated to a ``vector bundle'' over $C$, with
$s : C\rightarrow {\mathcal  E}_{{\mathcal  P},C}, \
s(c) = (c,c)$, as the zero section.

\begin{defin}
\label{Rene}
Let us define the ``Thom-space'' of the ``disc''-bundle
$p : {\mathcal  E}_{{\mathcal  P},C}\rightarrow C$
as the one-point compactification
$$
{\rm Thom}({\mathcal  E}_{{\mathcal  P},C}) :=
({\mathcal  E}_{{\mathcal  P},C}\setminus{\dot {\mathcal  E}}_{{\mathcal  P},C})
\cup \{\infty\} .
$$
\end{defin}
\begin{prop}
\label{Thom} Under assumptions {\bf A1} and {\bf A2},
the second projection map $\pi_2 :\Delta({\mathcal  P})^2\rightarrow
\Delta({\mathcal  P})$ induces a homeomorphism
$$
{\rm Thom}({\mathcal  E}_{{\mathcal  P},C}) \stackrel{\cong}{\longrightarrow}
\Delta({\mathcal  P})/\Delta({\mathcal  R}) \ .
$$
\end{prop}
{\bf Proof}: The projection $\pi_2$ induces a continuous map
$\alpha : ({\mathcal  E}_{{\mathcal  P},C}\setminus{\dot{\mathcal  E}}_{{\mathcal  P},C})
\rightarrow \Delta({\mathcal  P})\setminus\Delta({\mathcal  R})$
and a map $\overline{\alpha} : {\rm Thom}({\mathcal  E}_{{\mathcal  P},C})\rightarrow
(\Delta({\mathcal  P})\setminus\Delta({\mathcal  R}))\cup \{\infty\}$.
Note that $\Delta({\mathcal  R})$ is compact which
implies that the letter space is homeomorphic to
$\Delta({\mathcal  P})/\Delta({\mathcal  R})$. The map $\alpha$ is
\ 1--1 \ because $C$ is an antichain. It remains to be checked that the
inverse function $\alpha^{-1}$ is also continuous. This will follow
from the continuity of of the function $\beta = p\circ\alpha^{-1}
: \Delta({\mathcal  P})\setminus\Delta({\mathcal  R})\rightarrow C$.
Note that if $x = \lambda_1 x_1 +\ldots + \lambda_n x_n
\in\Delta({\mathcal  P})\setminus\Delta({\mathcal  R}), \lambda_i > 0$,
then $\beta(x)=y$ where $y$ is the unique
element in the intersection of the chain $Z := \{x_1<\ldots < x_n\}$
and antichain $C$. Let $U\subset C$ be a neighborhood of
$y$ in $C$. It follows from {\bf A2} that ${\mathcal  P}\setminus U$
is a closed subposet of ${\mathcal  P}$, hence $\Delta({\mathcal  P}\setminus U)$
is a closed subspace of $\Delta({\mathcal  P})$.
If $V:=\Delta({\mathcal  P})\setminus\Delta({\mathcal  P}\setminus U)$,
then $V$ is a neighborhood of $x$ and $\beta(V)\subset U$.
This means that $\beta$ is continuous which completes the proof
of the proposition. \hfill$\square$

\medskip
In the most interesting cases, $p : {\mathcal  E}_{{\mathcal  P},C}
\rightarrow C$ is indeed homeomorphic to the disc bundle
of an actual vector bundle over $C$.
In this case the ``Thom-space'' from definition~\ref{Rene}
can be described as the {\em twisted smash product} in the
sense of the following definition.

\begin{defin}
\label{twist1}
Suppose that ${\mathcal  E}$ is a fibre bundle pair over $B$ with total pair
$(E,\dot E)$, fiber pair $(F,\dot F)$ and a projection
$p : E\rightarrow B$, cf. \cite{Spanier}, p.256.
For example ${\mathcal  E}$ can be the disc and sphere bundle pair
$(D(V), S(V))$ associated to a vector bundle $V\stackrel{p}{\rightarrow}B$.
Suppose $\dot B\subset B$ is a closed subspace. Then
the twisted smash product
$$
(B,\dot B)\rtimes_{\mathcal  E} (F,\dot F)
$$
of pairs $(B,\dot B)$ and $(F,\dot F)$ relative ${\mathcal  E}$, is the
quotient space $E/(\dot E \cup p^{-1}(\dot B)$.
\end{defin}

\begin{cor}
\label{twist2}
Suppose that the compact bundle pair
$({\mathcal  E}_{{\mathcal  P},C}, \dot{\mathcal  E}_{{\mathcal  P},C})
\rightarrow C$ from definition~\ref{Rene} is isomorphic to the
disc and sphere bundle pair $(D(V),S(V))\rightarrow C$
associated to some $m$-dimensional vector bundle over $C$.
Assume that $V$ can be extended to a vector bundle
$W$ over a compact space $\widetilde C\supset C$. Then,
$$
\Delta({\mathcal  P})/\Delta({\mathcal  R}) \cong {\rm Thom}({\mathcal  E}_{{\mathcal  P},C})
\cong (\widetilde C, C)\rtimes_{\mathcal  W} (D^m,S^{m-1})
$$
where ${\mathcal  W} = (D(W),S(W))$ is the disc and sphere bundle
pair associated to $W$.
\end{cor}

\begin{cor}
\label{twist3}
If the bundle $V$ in the corollary~\ref{twist2} is trivial, then
$W$ can be taken to be the trivial bundle over
$\bar C := C\cup {\infty}$ in which case
$$
\Delta({\mathcal  P})/\Delta({\mathcal  R}) \cong {\rm Thom}({\mathcal  E}_{{\mathcal  P},C})
\cong (\bar C,\{\infty\}) \wedge (D^m,S^{m-1})
$$
where $"\wedge "$ is the usual smash product of pointed spaces.
\end{cor}
For the sake of completeness we formulate one more corollary of
proposition~\ref{Thom} which shows explicitly the connection
of this result with the equation~(\ref{wedge}).

\begin{cor}
Suppose that the topological poset ${\mathcal  P}$ satisfies,
besides {\bf A1} and {\bf A2}, the condition that
the antichain $C$ is finite. Then,
$$
\Delta({\mathcal  P})/\Delta({\mathcal  R})\cong {\rm Thom}({\mathcal  E}_{{\mathcal  P},C})
\cong \bigvee_{c\in C} \Sigma(\Delta({\mathcal  P}_{< c})\ast
\Delta({\mathcal  P}_{> c}))
$$
\end{cor}

\section{Applications}
\label{Appl}
\subsection{Grassmannian posets}
\label{Grass}
\begin{theo}
\label{Vas1}
(Vassilev, \cite{Vas91, Vas94})
Let $K = R, C, Q$ be one of the classical fields and
let  $G_k(K^n)$ be the Grassmann manifold of all linear,
$k$-dimensional subspaces of $K^n$.
Let ${\mathcal  G}_n(K) = (G(K^n),\subseteq)$
be the associated Grassmannian poset where  $G(K^n)$
is a disjoint sum of Grassmannians
$$ G(K^n) := \coprod_{i=0}^{n} G_i(K^n) $$
\noindent
Let $\widetilde{\mathcal  G}_n(K) = (\widetilde G(K^n),\subseteq)$ be the
truncated Grassmannian poset where $\widetilde G(K^n) :=
\coprod_{i=1}^{n-1} G_i(K^{n}) $.
Then,
$$
\Delta(\widetilde{\mathcal  G}_n(K)) \cong S^{\binom{n}{2}d+n-2}
\quad \mbox{\rm where} \ \ d = {\rm dim}_R(K).
$$
\end{theo}
\noindent
{\bf Proof:} We want to illustrate the technique developed
in this paper so we concentrate on the proof that
the order complex $\Delta(\widetilde{\mathcal  G}_n(K))$ has the
same homotopy type as the sphere of dimension $\binom{n}{2}d+n-2$.
The proof that  $\Delta(\widetilde{\mathcal  G}_n(K))$ is
actually a sphere can be completed along the lines of \cite{Vas91},
see also the section~\ref{Elem} for a direct, elementary proof.
To simplify the notation we prove the result in the case
$K = R$, in other two cases the proof is completely
analogous.  Let $Z\in G_1(R^n)$ and let
$C:={\rm Co}(Z)$ be the space of all complements of
$Z$ in ${\mathcal  G}_n(R)$. Obviously,
${\rm Co}(Z) = \{L\in G_{n-1}(R^n)\mid Z\not\subseteq L\}$.
It is easy to see that the space $C$ is homeomorphic to
$R^{n-1}$ so the one-point compactification $\bar C$ of
$C$ is homeomorphic to $S^{n-1}$.   Let ${\mathcal  G'}_n(R) =
\widetilde{\mathcal  G}_n(R)\setminus {\rm Co}(Z)$.
Then by proposition~\ref{4-3} and corollary~\ref{twist3},
$$ \Delta(\widetilde{\mathcal  G}_n(R)) \simeq
\Delta(\widetilde{\mathcal  G}_n(R))/\Delta({\mathcal  G'}_n(R)) \simeq
(\bar C,\infty) \wedge \Sigma(\bar{\mathcal  G}_{n-1}(R))
\simeq S^{n-1}\wedge \Sigma(\widetilde{\mathcal  G}_{n-1}(R))
$$
Obviously $\Delta(\widetilde{\mathcal  G}_2(R)) \simeq S^1$ so
an induction based on the homotopy recurrence relation above
yields the desired formula
$$
\Delta(\widetilde{\mathcal  G}_n(R)) \simeq S^{\binom{n}{2}+n-2} \ .
$$

\medskip
The idea of the proof of theorem~\ref{Vas1} is quite general
and can be applied in many other situations.
Here is another example where we compute the homotopy
type of an oriented Grassmannian poset.

\begin{theo}
Let  $G_k^{\pm }(R^n)$ be the Grassmann manifold of all linear,
oriented, $k$-dimensional subspaces of $R^n$.
Let ${\mathcal  G}_n^{\pm}(R) = (G^{\pm}(R^n),\subseteq)$
be the associated Grassmannian poset
and $\widetilde{\mathcal  G}_n^{\pm}(R) = (\widetilde G^{\pm}(R^n),\subseteq)$ the
associated truncated poset where $\widetilde G^{\pm}(R^n) :=
\coprod_{i=1}^{n-1} G_i^{\pm}(R^{n}) $ and
$G^{\pm}(R^n) := \widetilde G^{\pm}(R^n)\cup G_0^{\pm}(R^n)
\cup G_n^{\pm}(R^n)$.
Then the homotopy type of the truncated
poset is equal to the wedge of $2^{n-2}$
spheres of dimension $\binom{n}{2} + n -2$,
$$
\Delta(\widetilde{\mathcal  G}_n^{\pm}(R)) \simeq
\bigvee_{j=1}^{2^{n-2}} S_j^{\binom{n}{2}+n-2}.
$$
\end{theo}
\noindent
{\bf Proof:} The proof is similar to the proof of theorem~\ref{Vas1}.
Choose $Z\in G_1^{\pm}(R^n)$ and denote by $C = {\rm Co}(Z)$ the space of all
complements of $Z$ in ${\mathcal  G}_n^{\pm}(R)$.
This space is seen to be a disjoint union of `positive' and
`negative' hyperplanes
$${\rm Co}(Z) = \{L\in G_{n-1}^{\pm}(R^n) \mid Z\not\subseteq
L\} \cong R^{n-1} \sqcup R^{n-1}.$$
\noindent
>From here we deduce the recurrence formula~(\ref{grasgras_lat})
in section~\ref{main} and the desired formula follows by
induction. \hfill $\square$

\medskip
\noindent
Let $I := \{i_{\nu}\}_{\nu = 1}^k \subseteq [n]$ be a
subset of $[n] = \{0,1,\ldots ,n\}$.
Define
$$
{\mathcal  G}_I(K) := \coprod_{\nu =1}^k \ G_{i_\nu}(K^n)
$$
as the ``rank selected'' subposet of ${\mathcal  G}_n(K)$.
A natural problem is to express the homotopy type
of $\Delta({\mathcal  G}_I(K))$ in terms of ``standard''
spaces and constructions. Here is an example.

\begin{prop}
Let $I = \langle k\rangle = \{1,\ldots ,k\}$.
Let ${\mathcal  G}_I(R) = \coprod_{\nu =1}^k~G_{\nu}(R^n)$
be the rank selected subposets of ${\mathcal  G}_n(R)$
associated to $I$. Then,  cf. definition~\ref{twist2},
$$
\Delta({\mathcal  G}_I(R))\cong (G_k(R^n),G_{k-1}(R^{n-1}))
\rtimes_{\mathcal  W} \ (D^m, S^{m-1})
$$
where $m = \binom{k}{2} + k -1$ and
${\mathcal  W}$ is the disc bundle which arises if each
fibre $L$ in the tautological $k$-plane bundle
over $G_k(R^n)$ is replaced by the cone over
$\Delta(\widetilde{\mathcal  G}_k(R)) \cong D^m$.
\end{prop}

\begin{rem}
In the case of general rank selected posets ${\mathcal  P} = {\mathcal  G}_I(K)$,
given $Z\in {\mathcal  P}$, the space ${\rm Co}(Z)$ of complements
of $Z$ is not necessarily an antichain.
The answer to the question {\bf Q3} from section~\ref{Homo}
gets more complicated, compare \cite{Bjo89-1},
and may require new ideas.
\end{rem}

\subsection{Configuration posets}
\label{Conf}

Configuration posets ${\rm exp}_n(X)$ and configuration spaces
$F(X,n)$ and $B(X,n)$, of labeled and unlabeled
points in the space $X$ respectively, were  defined in example~\ref{exam2}.
The following construction has been introduced
by Vassiliev under the name {\em geometric resolution} of
configuration spaces ${\rm exp}_n(M)$ and
$B(M,r)$.
Suppose that a manifold or more generally a finite CW-complex $M$
is generically embedded in the space $R^N$ of very large dimension
$N$. Let ${\rm Conv}_r(M)$ be the union of all
(closed) $(r-1)$-dimensional simplices with vertices
in the embedded space $M$. The genericity of the embedding
means that two simplices spanned by different sets of
vertices must have disjoint interiors.
It is easy to observe that
${\rm Conv}_n(M) \cong \Delta({\rm exp}_n(M))$, more precisely
the order complex  $\Delta({\rm exp}_n(M))$ can be seen as
the barycentric subdivision of ${\rm Conv}_n(M)$.
This is precisely the context in which order complexes of
continuous posets arose in Vassiliev's work.
For example he proved in \cite{Vas97} the following
theorem.

\begin{theo}
\label{Vas2}
The space ${\rm Conv}_n(S^1) \cong \Delta({\rm exp}_n(S^1))$
is homeomorphic to $S^{2n-1}$.
\end{theo}
Vassiliev actually proved a little more by showing that spaces above are
$PL$-ho\-meo\-mor\-phic.
Motivated by this example, he asked in his lecture
at the workshop ``Geometric Combinatorics'',
MSRI, February 1997, whether a similar formula holds for other spaces
(manifolds) $M$. In other words is it always true that
$$
\Delta({\rm exp}_n(X)) \cong X\ast \ldots\ast X = X^{\ast(n)}
$$
where $X^{\ast(n)}$ denotes the $n$-fold join of the space
$X$. Here we show, theorem~\ref{conf1}, corollary~\ref{conf2} and
proposition~\ref{conf3},
that the homology type of the space
$\Delta({\rm exp}_n(M))$ is very closely related to the
homology type of the space $B(M,n)$. As a consequence we obtain
a proof of theorem~\ref{Vas2} and show that already in the
case $X = S^2, n \geq 2$ the answer to the question raised
by Vassiliev is negative.

We start with a proof of $\Delta({\rm exp}_n(S^1))\cong S^{2n-1}$.
Vassilev gave in \cite{Vas97} an elegant and short proof of this
result. Another elementary proof is given in section~\ref{Elem}
while the following proof illustrates in the first place
some of the tools developed in
this paper.

\medskip
\noindent
{\bf Proof of theorem~\ref{Vas2}:} It is not difficult to check that
$\Delta({\rm exp}_n(S^1))$ is a $PL$-manifold so,
in light of the fact that the Poincar\' e conjecture holds for
$PL$--manifolds of dimension $n\geq 5$
(the case $n=2$ is established by a direct argument),
we concentrate on the proof that
$\Delta({\rm exp}_n(S^1))$ has the correct homotopy type.
Let $x_0\in S^1 = B(S^1,1)\subset {\rm exp}_n(S^1)$.
The set ${\rm Co}(\{x_0\})$ of all complements of
$\{x_0\}$ in the lattice ${\rm exp}_n(S^1) \cup \{\hat{0}, \hat{1}\}$
is $B(S^1\setminus\{x_0\}, n)$. The space $B(S^1\setminus\{x_0\}, n)
\cong B(R^1, n)$ is an open, $n$-dimensional convex set so
its one--point compactification is homeomorphic to $S^n$.
Finally, by the results of sections~\ref{1and1} and \ref{2and3},
$$
\Delta({\rm exp}_n(S^1)) \simeq S^{n}\wedge \Delta(\widetilde{\mathcal  B}_{n-1})
\simeq S^{2n-1}
$$
and the result follows. \hfill$\square$

\medskip
The idea of the proof of theorem~\ref{Vas2} can be obviously used for more
general configuration posets.

\begin{defin}
\label{bundle}
The configuration space $F(Y,n)$ of all labeled $n$-element
subsets of $Y$ is a free $S_n$-space. Let $V$ be the standard
representation of $S_n$ i.e. the $(n-1)$-dimensional
representation obtained by subtracting the trivial $1$-dimensional
representation of $S_n$ from the permutation representation.
Let $\xi$ be the vector bundle
\[
R^{n-1}\longrightarrow F(Y,n)\times_{S_n} V\longrightarrow B(Y,n)
\]
and let us define its Thom-space as its one-point compactification
$$
{\rm Thom}_n(Y) := (F(Y,n)\times_{S_n} V)\cup \{\infty\} .
$$
\end{defin}

\begin{defin}
\label{admissible}
Let ${\mathcal  P} = {\rm exp}_n(X)$ be the $n^{{\rm th}}$ configuration
poset of $X$ and $x_0\in X$. Let $C := {\rm Co}(\{x_0\})$
be the space of all complements of $\{x_0\}$ in the lattice
$\hat{\mathcal  P} = {\mathcal  P}\cup\{\hat 0, \hat 1\}$ and
${\mathcal  R}:= {\mathcal  P}\setminus C$.
The pointed space $(X,x_0)$ is called {\em admissible} if the
inclusion $\Delta({\mathcal  R}) \hookrightarrow \Delta({\mathcal  P})$
is a cofibration, cf. section~\ref{2and3}. For example
every pair $(M, p)$ is admissible if $M$ is a smooth,
compact manifold and $p\in M$.
\end{defin}

\begin{theo}
\label{conf1}
Suppose that $(X, x_0)$ is a compact, admissible space.
Then
$$
\Delta({\rm exp}_n(X)) \simeq {\rm Thom}_n(X\setminus\{x_0\})
$$
\end{theo}

\medskip
\noindent
{\bf Proof}: The compactness of  $\Delta({\rm exp}_n(X))$ follows
from the compactness of $X$. The space $C = {\rm Co}(\{x_0\})
\cong B(X\setminus\{x_0\}, n)$ is an open antichain in
$\Delta({\rm exp}_n(X))$ so both conditions $A1$ and $A2$
of proposition~\ref{Thom} are satisfied. The bundles
${\mathcal  E}_{{\mathcal  P},C}$ and $\dot{\mathcal  E}_{{\mathcal  P},C}$
from definition~\ref{Rene} are identified as the disc and the sphere
bundle associated to the bundle from the definition~\ref{bundle}.
The result is therefore a consequence of proposition~\ref{Thom}.
\hfill$\square$

\begin{cor}
\label{conf2}
The order complex $\Delta({\rm exp}_n(S^m))$, has the
homotopy type of the space
$$
{\rm Thom}_n(R^m) = (F(R^m,n)\rtimes_{S_n} V)\cup \{\infty\}
$$
By the Poincar\' e-Lefschetz duality, taking the (co)homology
with $Z_2$-coefficients,
$$
\widetilde{H}_p(\Delta({\rm exp}_n(S^m))) \cong H^{(m+1)n-p-1}(B(R^m,n)).
$$
\end{cor}

\begin{prop}
\label{conf3}
$\Delta({\rm exp}_n(S^2))$ does not have the homotopy type of a sphere
for $n\geq 2$.
\end{prop}
\noindent
{\bf Proof:} By corollary~\ref{conf2}
$$H_p(\Delta({\rm exp}_n(S^2))) \simeq H^{3n-p-1}(B(R^2,n)).$$

The cohomology of the configuration space $B(R^2,n)$ with
$Z_2$ coefficients has been computed by D.B.~Fuchs,
\cite{Fuchs, Vas94}. He
proved that the dimension of $H^k(B(R^2,n), Z_2)$ over $Z_2$
is equal to the number of representations
of the form $n = 2^{\alpha_1}+\ldots + 2^{\alpha_{n-k}},
\, \alpha_i\geq 0$, where two representations that differ
only by the order of summation are considered to be equal.
>From here we deduce that $H^0(B(R^2,n), Z_2) \cong
H^1(B(R^2,n), Z_2) \cong Z_2$, hence the reduced homology
of $\Delta({\rm exp}_n(S^2))$ is nontrivial in two dimensions
i.e. this space cannot be homotopic to a sphere.
\hfill$\square$

\medskip
The following simple corollary reveals the connection of theorem~\ref{conf1}
and its consequences with the interesting problem of $n$-neighborly
submanifolds of $R^N$, \cite{Vas97c}. Note that in this statement,
the embedding of $S^2$ is not necessarily smooth, let alone stably
$n$-supported in the sense of \cite{Vas97c}. This suggests that
a more general problem of finding $n$-neighborly embeddings
of more general spaces, say simplicial complexes, may be interesting.

\begin{cor} Suppose that the sphere $S^2$ is topologically
embedded in $R^N$ in such a way that for any collection $\Sigma$ of $r$
different points on the sphere, there exists a hyperplane in $R^N$ such
that $\Sigma\subset H$ and $S^2\setminus \Sigma$ is in an open
halfspace determined by $H$. Then  $N\geq 3n$.
\end{cor}

\medskip
\noindent
{\bf Proof:} Every embedding of the sphere $S^2$ in $R^N$ which is
$n$-neighborly in the sense above, leads to an embedding of the order
complex ${\rm exp}_n(S^2)$ into the same space. By the proof
of the corollary~\ref{conf3} we know that $H_{3n-1}({\rm exp}_n(S^2))
\neq 0$, hence $3n\leq N$.

\subsection{Elementary proofs}
\label{Elem}
In this section we present direct, elementary proofs of
two motivating results of Vassiliev, theorems \ref{Vas1} and \ref{Vas2}
of this paper. I am obliged to W.~Thurston for the idea of the first,
to G.~Kalai for the idea of the second proof and to
B.~Shapiro for the information related to convex curves in $R^n$.

\medskip
\noindent
{\bf Second proof of theorem~\ref{Vas1}:}
Let ${\rm Symm}(n)$ be the vector space of all
symmetric, $n\times n$ matrices with real entries.
The dimension of this space is $\binom{n}{2} + n$.
Given a matrix $A := [a_{ij}]$ let
$\lambda_1\leq \lambda_2\leq\ldots\leq\lambda_n$
be the ordered sequence of its eigenvalues and let
$u_1, u_2,\ldots , u_n$ be a sequence of linearly
independent vectors such that $u_i$
is an eigenvector associated to the eigenvalue $\lambda_i$.
Let $L_1\subset L_2\subset\ldots\subset L_n$
be the flag associated to this sequence,
$L_i := {\rm span}\{u_j \mid j = 1,\ldots ,i\}$.
Note that this flag is not well defined if some of the eigenvectors
coincide. On the other hand if
$M :=  {\rm Symm}(n)\setminus \{\alpha I\mid \alpha\in R\}$
is the space of all matrices which are not multiples of the unit matrix,
then the map
$$
\Phi : M  \longrightarrow G_1(R^n) \ast \ldots \ast G_{n-1}(R^n)
$$
\noindent
defined by
$$
A := [a_{ij}]  \mapsto \frac{\lambda_2 - \lambda_1} {\lambda_n - \lambda_1}
L_1 +  \frac{\lambda_3 - \lambda_2} {\lambda_n - \lambda_1} L_2 +\ldots +
\frac{\lambda_n - \lambda_{n-1}} {\lambda_n - \lambda_1}
L_{n-1}
$$
is well defined. This map is not one--to--one. More precisely,
$\Phi(A) = \Phi(B)$ iff there exist $\alpha > 0$ and $\beta\in R$
such that $B = \alpha A + \beta I$ or in other words
$\Phi$ is constant on orbits of the group
$R^+\times R$ which acts on $M$
by the formula $(\alpha,\beta)A := \alpha A + \beta I$.
We conclude that the image of this map,
${\rm Image}(\Phi) = \Delta(\widetilde{\mathcal  G}_n(R))$,
is diffeomorphic  to the orbit manifold $M/(R^+\times R)$.
The orbits are easily identified as `vertical' two dimensional
halfplanes with the origin removed. Hence,
if $H$ is the codimension one linear subspace of
${\rm Symm}(n)$  orthogonal to $I$, then
the $(\binom{n}{2} + n -2)$-dimensional
unit sphere $S(H)$ in $H$ intersects each orbit in exactly one point.
This completes the proof.

\bigskip

\noindent
{\bf Second proof of theorem~\ref{Vas2}:}
This proof relies on some characteristic properties of
so called {\em closed convex curves} in $R^{2n}$, \cite{Sch54,
SedSha}.
We outline the main idea of the proof leaving details to
the reader.
A simple, closed curve $\gamma : S^1 \rightarrow R^{2n}$
is called {\em convex} if the total multiplicity
of its intersection with any affine hyperplane does
not exceed $2n$. An example is the Carath\' eodory curve
$$
\Gamma : t \mapsto (\sin t, \cos t,
\sin (2t), \cos (2t), \ldots , \sin (nt), \cos (nt)).
$$
Let $C := {\rm conv}({\rm Im}(\Gamma))$ be the convex hull
of this curve and let $D := \partial (C)\cong S^{2n-1}$.
All we have to show is that every point in
$D$ can be expressed uniquely as a convex combination
of the form $\alpha_1 x_1 + \ldots + \alpha_n x_n$ where
$x_1, x_2, \ldots , x_n$ are distinct points in
${\rm Im}(\Gamma)$ and conversely that every such convex
combination determines a point in $D$.
Both facts follow from the following observation
of Schoenberg, \cite{Sch54}, who proved that
the convex hull of every closed convex curve $\gamma$ is the
intersection of a family of closed halfspaces determined to support
hyperplanes which are tangent to   $\gamma$ at $n$ distinct
points. Note that, if $C$ is seen as a continuous
analogue of the {\em cyclic} polytope,
then this result of Schoenberg can be related to the
well known Gale's evenness condition, \cite{MullShe, Zie95}.
Indeed, hyperplanes tangent to $\gamma$ at $n$ distinct points
are spanned by $n$ pairs of ``infinitesimally closed'' points
on this curve.

\end{document}